\documentclass[12pt]{article}
\usepackage{graphicx}
\usepackage{amsmath}
\usepackage{amssymb}
\usepackage{url} 
\usepackage{amsthm}
\usepackage{tikz-cd}
\usepackage{hyperref}
\usepackage{quickmath}
\usepackage{epsfig}
\usepackage{latexsym}      

\usepackage{lineno}
\usepackage{changes}

\newtheoremstyle{letter}%
  {}
  {}
  {}
  {}
  {\bfseries}
  {.}
  { }
  {}%

\newcommand{\cah}{\mathcal{H}}

\newcommand{\bbn}{\mathbb{N}}

\theoremstyle{letter}
\newtheorem{thmx}{Theorem}

\newenvironment{myproof}[1][\proofname]{%
  \proof[ \textbf{Proof}]%
}{\endproof}

\newtheorem{thmz}{Theorem}
\date{}
\title{Rigidity Properties for Hyperbolic Generalizations}
\author{Brendan Burns Healy}
\begin{document}

\maketitle
\abstract{
We make a few observations on the absence of geometric and topological rigidity for acylindrically hyperbolic and relatively hyperbolic groups.  In particular, we demonstrate the lack of a well-defined limit set for acylindrical actions on hyperbolic spaces, even under the assumption of universality.  We also prove a statement about relatively hyperbolic groups inspired by a remark by Groves, Manning, and Sisto in \cite{sisgroman} about the quasi-isometry type of combinatorial cusps.  Finally, we summarize these results in a table in order to assert a meta-statement about the decay of metric rigidity as the conditions on actions on hyperbolic spaces are loosened.
}


\begin{section}{Introduction}

Gromov-hyperbolic spaces are a core focus in geometric group theory, in part because of how they behave under deformation.  This is a stark contrast from generic metric spaces and their large-scale properties.  For example, the geodesic ray definition for boundaries of spaces is not always homeomorphically rigid under quasi-isometry, famously failed by $\CAT(0)$ spaces as shown by Croke and Kleiner in \cite{croke}. Therefore the visual boundary of a $\CAT(0)$ group is not well-defined, however we find that it will be if the space $X$ is hyperbolic.  Specifically, any quasi-isometry $X \rightarrow Y$ induces a homeomorphism $\partial X \cong \partial Y$.  Furthermore, then, if a group acts properly, cocompactly, and by isometries on a hyperbolic space $X$, by the \v Svarc-Milnor lemma, we can even make sense of $\partial G$.

In an effort to make more broad statements, geometric group theoriests often loosen the requirements of a `geometric' action.  In particular we may consider groups which act acylindrically and nonelementarily on hyperbolic spaces, a class which are aptly named \emph{acylindrically hyperbolic}.  Much of the known machinery and results regarding this class is available in \cite{osin}.  We also consider a class of groups introduced by Gromov with important early attention given by Bowditch and Farb  - \emph{relatively hyperbolic} groups.  Specifically these are groups that act in a `geometrically finite' way on hyperbolic spaces that may be thought of as Cayley graphs with negatively curved cusps added.  The geometrical finiteness in this action refers to the fact that although the quotient space is no longer compact, it has a finite number of \emph{ends}, in the sense of \cite{geog}.

Although we can define what the boundary of the group should mean in this latter case, we find that extending that to acylindrically hyperbolic groups, which is a strict generalization of relatively hyperbolic groups, is problematic.  We demonstrate that even under some additional assumptions the limit sets of these actions do not achieve a consistent shape.  

This is stated as follows, with definitions of terms coming in section 2.

\begin{thmx}
There exist acylindrically hyperbolic groups $G$, which admit two different universal actions $G \curvearrowright X_i$, such that in the representations 
\[ \rho_1: G \rightarrow \emph{Isom}(X_1), \ \rho_2: G \rightarrow \emph{Isom}(X_2) \]
the limit sets $\Lambda_1(G)$ and $\Lambda_2(G)$ are not homeomorphic.

\end{thmx}

In fact, the actions will be by any closed surface group, both on the space $\H^3$, where one copy will be identified with the universal cover of a 3-manifold, in which the surface sits as a normal subgroup and has full limit set, and the other will have the action induced by a geodesically embedded copy of $\H^2$ with the action extending the natural one by deck transformations, with limit set a circle.

Recent work of Abbott, Balasubramanya, and Osin in \cite{ABO}, has generated the idea of a \emph{largest} such action, though this is only possible in the class of cobounded actions. We note however, that these actions are not guaranteed to exist.  This is because universal actions themselves are not always guaranteed to exist by \cite{abbott}.

One advantage of working in the setting of geometric actions on hyperbolic spaces is the quasi-isometry invariance of the acted-on spaces and their boundaries.  We can recover some well-defined notion of boundary for geometrically finite actions on hyperbolic spaces, so we may ask if we retain the quasi-isometry invariance as well. We should expect that in the `well-behaved' portion of the hyperbolic cusped space, which is the preimage of a compact portion of the quotient, chosen so that the lifts of the cusps don't intersect, we do see a nice invariance.  So it must be then that if quasi-isometry invariance breaks down, it is in the cusps.  It is a fact asserted in \cite{sisgroman} that one can choose the shape of these cusps carefully (or carelessly, depending on your viewpoint) to force them away from being in the same QI class.  We prove this fact rigorously, to obtain a statement about geometrically finite actions.

\begin{thmx} 
Any relatively hyperbolic group with infinite peripheral subgroups acts as in \cite{bowditch} on hyperbolic spaces that are not equivariantly quasi-isometric.
\label{thmb}
\end{thmx}

\vspace{.1in}

The author would like to thank Genevieve Walsh, Robert Kropholler, and Daniel Groves for many helpful conversations, as well as the reviewer for helpful suggestions.

\end{section}

\begin{section}{Universal Acylindrical Actions}
To define acylindrically hyperbolic groups, we define what it means for an action to be acylindrical. 

\defi
 An metric space action $G \curvearrowright X$ is called \emph{acylindrical} if for every $\epsilon > 0$ there exist $R(\epsilon ),N(\epsilon ) > 0$ such that for any two points $x,y \in X$ such that $d(x,y) \geq R$, the set
$$ \{ g \in G \mid d(x,g.x) \leq \epsilon, d(y,g.y) \leq \epsilon \}$$
has cardinality less than $N$.
\edefi

We need more hypotheses on a group than simply acting acylindrically on a hyperbolic space, however, as all groups admit such an action.  The trivial action of any group on a point, a hyperbolic space, is acylindrical.

\defi
A group $G$ is called \emph{acylindrically hyperbolic} if it admits an acylindrical action on a hyperbolic space which is not elementary; that is, it has a limit set inside the boundary of the space of cardinality strictly greater than 2.
\edefi

In restricting to this class, we obtain a more interesting class of groups.  Indeed, we omit some groups we feel in some natural sense, shouldn't be negatively curved.  A quick fact available in \cite{osin}, for example, tells us that any group that decomposes into the direct product of two infinite factors is not acylindrically hyperbolic.  For some subclasses,such as right angled Artin groups (see \cite[Section~8]{osin}, \cite{sisto} and \cite{capsag}), this, together with being virtually cyclic, is a complete obstruction to acylindrical hyperbolicity.
We would like to know if, similar to hyperbolic and relatively hyperbolic groups, these groups admit some well-defined notion of boundary or limit set.  This question is also being studied by Abbott, Osin, and Balasubramanya, who in \cite{ABO} develop what they term a \emph{largest} action, which is necessarily also cobounded.  Although here we will not look at cobounded actions, we do use one of their conditions, which is that our actions will be \emph{universal}.

\defi
Let $G$ be an acylindrically hyperbolic group.  An element $g \in G$ is called a \emph{generalized loxodromic} if there's an acylindrical action $G \curvearrowright X$ for $X$ hyperbolic such that $g$ acts as a loxodromic.  An individual isometry $g$ is a \emph{loxodromic} if for some basepoint $x_0 \in X$, the map $\Z \rightarrow X$ defined by $n \mapsto g^n . x_0$ is a quasi-isometry.
\edefi

\defi
For an acylindrically hyperbolic group, an action $G \curvearrowright X$ is called \emph{universal} if it is acylindrical, $X$ is hyperbolic, and all generalized loxodromics act as loxodromics.
\edefi

Universal actions are a natural setting to consider our question, as we can easily change a given action if we force a generalized loxodromic to act elliptically.  Even with universality, however, we do not get a well-defined boundary.  First we note that a subgroup of an acylindrically hyperbolic group will inherit that property if the induced sub-action remains non-elementary.

We will also need to note that geometric actions are acylindrical. This is observed in \cite{osin}, but we provide a proof here for completeness.

\lem \cite{osin}
If a group action is geometric then it is also acylindrical.
\elem

\pf
Suppose $G \curvearrowright X$ geometrically.  Let $K \subset X$ be a compact fundamental domain for this action.  Set $d = diam(K)$.  We note by cocompactness that for any $x,y \in X$, there exists a group element $h \in G$ such that 
\[ d(x, h.y) \leq d . \]
We make one more claim, that is due to the action being by isometries.  We claim that for all $\epsilon > 0, y \in X$ and $h$ as above,
\[ \{g \in G | d(y,g.y) \leq \epsilon \} = \{g \in G | d(h.y, h.(g.y) \leq \epsilon \}. \]

Now, for $\epsilon > 0$, pick $R(\epsilon) > d$.  For any two points $x,y$, we can choose $g$ such that $d(x,g.y) \leq d$, i.e. that both $x,g.y$ belong to the same translate of $K$.  Without loss of generality, assume this translate is $K$ itself.  Then the set 
\[ \{ g \in G | d(x,g.x) \leq \epsilon, d(y,g.y) \leq \epsilon \} \]
is exactly equal to the set 
\[ \{ g \in G | d(x,g.x) \leq \epsilon, d(h.y,h.(g.y)) \leq \epsilon \}. \]
This set is a subset of the set of elements which translate $K$ to a tile at distance a maximum of $\epsilon$ away, which is bounded because the group action is proper.  This bound is a function of $\epsilon$, so let this bound serve as $N(\epsilon)$
\epf

The groups that we invoke for our non-uniqueness claim will be hyperbolic surface groups which will act on $\H^3$.  Accordingly, we need one more lemma, to do with hyperbolic geometry.

\lem
Let $\Gamma$ be a torsion-free Fuchsian group acting geometrically on $\H^2$.  Then for the natural isometric embedding of $\H^2 \hookrightarrow \H^3$, the induced action of $\Gamma \curvearrowright \H^3$ that comes from the inclusion $PSL(2,\R) \subset PSL(2,\C)$ is acylindrical.
\elem

\pf
Label by $X$ the original embedded copy of $\H^2$.  Let $x \in X$ be a point in this subspace.  Then $x$ belongs to some fundamental domain $K$ of the action $\Gamma \curvearrowright X$.  \\

\begin{figure}[h]
\begin{center}
\begin{tikzpicture}[scale=0.4]
\draw [thick] circle [radius=8];
\draw [out=330, in =210] (-8,0) to (8,0);
\draw [dashed, out=30, in =150] (-8,0) to (8,0);
\draw[out=310,in=205] (-1,1) to (1,1);
\draw [out=230, in=170] (1,1) to (2,-.2);
\draw [out=175, in= 60] (2,-.2) to (0,-1.1);
\draw [out=90, in=348] (0,-1.1) to (-1.3,-.2);
\draw [out=30, in =280] (-1.3,-.2) to (-1,1);

\node[color=red] at (2,-.2) {\textbullet};
\node[color=red] at (1,1)  {\textbullet};
\node[color=red] at  (-1,1) {\textbullet};
\node[color=red] at  (0,-1.1)  {\textbullet};
\node[color=red] at (-1.3,-.2) {\textbullet};

\draw [color= red, out=90, in=230] (2,-.2) to (4.8,6.4);
\draw [color= red, out=90, in=230, dashed] (1,1) to (3.25,7.31);
\draw [color= red, out=90, in=330, dashed] (-1,1)  to (-3.25,7.31);
\draw [color= red, out=90, in=315] (0,-1.1) to (-1.8,7.8);
\draw [color= red, out=90, in=330] (-1.3,-.2) to  (-4.8,6.4);

\draw [color= red, out=270, in=130] (2,-.2) to (4.8,-6.4);
\draw [color= red, out=270, in=130, dashed] (1,1) to (3.25,-7.31);
\draw [color= red, out=270, in=30, dashed] (-1,1)  to (-3.25,-7.31);
\draw [color= red, out=270, in=30] (0,-1.1) to (-1.5,-7.8);
\draw [color= red, out=270, in=30] (-1.3,-.2) to (-4.8,-6.4);
\end{tikzpicture}
\end{center}
\caption[Fundamental Domain for a Certain Non-cocompact Action on $\H^3$]{The convex hull of these geodesics serve as a fundamental domain}
\label{fixptlem}
\end{figure}
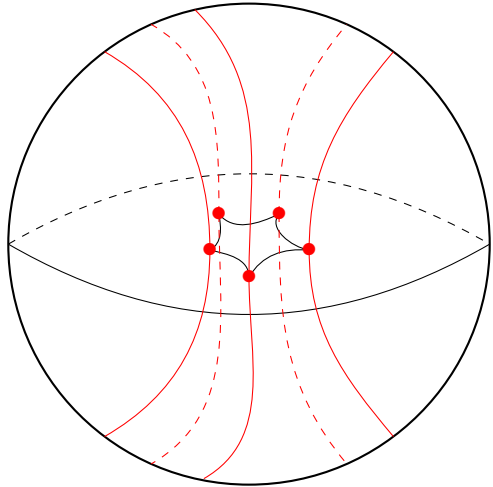

For any $\epsilon > 0$, let $N(\epsilon)$ be the (necessarily finite) number of translates of $K$ that intersect $\mathcal{N}_\epsilon(K)$.  That is to say $ N(\epsilon) = |S| $ where 
\[ S = \{ g \in \Gamma \ | \ gK \cap \mathcal{N}_\epsilon(K) \neq \emptyset \} \]
\newline
For any $g \in \Gamma$, which are all acting as loxodromics because $\Gamma$ is assumed to be torsion-free, because $X$ is the only totally geodesic copy of $\H^2$ that $\Gamma$ acts on geometrically, the geodesic axis lies entirely within $X$.  What this tells us is that for any point $z \in \H^3$, there exists a point $x \in X$, such that 
\[ d(z,g.z) \geq d(x, g.x) . \]  From this we can determine that 
\[ \{ g \in \Gamma | d(z,g.z) \leq \epsilon \} \subset \{ g \in \Gamma | d(x,g.x) \leq \epsilon \}. \]
\newline
The size of this right-hand set is bounded by $N(\epsilon)$, which thus also bounds the size of the left-hand set.  Because for any points $w,z \in \H^3$
\[ \{ g \in \Gamma | d(z,g.z) \leq \epsilon, d(w,g.w) \leq \epsilon \} \subset \{ g \in \Gamma | d(z,g.z) \leq \epsilon \} \]
we get that the action is acylindrical with constants $N(\epsilon)$ as chosen before, and any value of $R(\epsilon) > 0$.
\qedhere
\epf

We have the pieces now to state the following.

\begin{thmz}
There exist acylindrically hyperbolic groups $G$, which admit two different universal actions $G \curvearrowright X$, such that in the representations 
\[ \rho_1: G \rightarrow \emph{Isom}(X), \ \rho_2: G \rightarrow \emph{Isom}(X) \]
the limit sets $\Lambda_1(G)$ and $\Lambda_2(G)$ are not homeomorphic.
\end{thmz}

\pf
The space in question will be $\H^3$ and the group a closed surface group.

The following argument will work for the fundamental group of any closed surface of genus $\geq 2$.  However to be explicit, we will consider $G = \pi_1(\Sigma_2)$ where $\Sigma_2$ is a closed surface of genus 2.

Now, consider the action $G \curvearrowright \H^2$.  This action that of deck transformations, recognizing $\H^2$ as the universal cover, $\widetilde{\Sigma_2}$.  Because the quotient of this space is a closed manifold, the action is geometric, meaning it is acylindrical.  Furthermore, it has full limit set; that is to say $\partial G =  \partial \H^2 \cong S^1$.  Finally, every nontrivial element in this group action acts as a loxodromic, meaning it is a universal action.

By the lemma above, this action extends to an acylindrical action on $\H^3$, that has limit set $\Lambda(G)\cong S^1$, with all nontrivial elements continuing to act loxodromically.

Now we want to exhibit another universal action by this group on $\H^3$ with distinct limit set.  Let $\phi$ be a Pseudo-Anosov element of $MCG(\Sigma_2)$.  We can construct a hyperbolic 3-manifold, the geometry of which is given to us by \cite{thurston}, by taking the space $\Sigma_2 \times [0,1]$, and identifiying $\Sigma_2 \times \{0\}$ with $\phi(\Sigma_2) \times \{1\}$.  Denote this manifold by $M$.  We get a decomposition of $\pi_1(M) = \pi_1(\Sigma_2) \rtimes_{\phi^*} \Z$, where $\phi^* \in Aut(\pi_1(\Sigma_2))$ is induced by $\phi$.

Again because the quotient is a closed manifold, the natural covering space action $\pi_1(M) \curvearrowright \H^3$ is geometric, and therefore acylindrical.  Also by the geometric nature of the action, we get $\partial \pi_1(M) = \partial \H^3 \cong S^2$. 

We use a fact proved by Thurston in \cite{thurston}, Corollary 8.1.3 in Chapter 8, to assert that in fact the $\pi_1(\Sigma_2)$ has the same limit set as the entire group, by normality.  Specifically, $\Lambda({\pi_1(\Sigma_2)}) = \partial \H^3 \cong S^2$.

Now we need to know that all elements act as loxodromics.  Because the action is geometric (and acylindrical), none will act as parabolics.  Therefore, we need to rule out the possibility of elements acting elliptically.  However, because $\H^3$ is $\CAT(0)$, we know that any element acting elliptically will have a fixed point on the interior of $\H^3$.  We note that all nontrivial elements of $\pi_1(M) = \pi_1(\Sigma_2) \rtimes_{\phi^*} \Z$ are infinite order.  This implies that none can act elliptically.  If a nontrivial $g \in \pi_1(M)$ was elliptic, then it would fix a point, giving us an infinite number of elements, the powers of g, fixing a point, which violates the properness assumption of our action.
\newpage
Thus the induced action $G = \pi_1(\Sigma_2) \curvearrowright \H^3$ has the following properties:
\begin{itemize}
\item It is acylindrical
\item The space is hyperbolic
\item It has limit set $S^2$
\item All nontrivial elements act as loxodromics
\end{itemize}

Therefore this is a universal action with a distinct (homeomorphism type) of limit set for the group $G$.
\epf

\end{section}

\begin{section}{Geometrically Finite Actions}
There are many exisiting definitions in the literature for what it means for a group/subgroup combination to be a relatively hyperbolic pair.  A thorough review of these conditions and their equivalence is available in \cite{groman}.  For our purposes, we will use the following.

\defi
A group action $G \curvearrowright M$, for $M$ a compact metrizable space, is called a \emph{convergence action} if the induced action on the set of distinct triples 
$$ \{ (m_1, m_2, m_3) | m_i \in M, m_i \neq m_j \text{ for } i \neq j \}$$
is properly discontinuous.
\edefi

To this end, we want actions which are a certain kind of convergence action.

\defi
A convergence action is called \emph{geometrically finite} if every $m \in M$ is such that one of the following is true:
\begin{itemize}
\item The point $m$ is a bounded parabolic point, meaning it has infinite stabilizer acting cocompactly on $M \setminus \{m\}$.
\item The point $m$ is a conical limit point, meaning there exists a sequence $g_i, i \in \N$ of group elements, and distinct points $a,b \in M$ such that $g_i m \rightarrow a$ and $g_i m' \rightarrow b$ for all $m' \neq m$.
\end{itemize}
\edefi

We call a group acting on a hyperbolic space a \emph{geometrically finite} action, if its induced action on the boundary of that space is a geometrically finite convergence action.

\begin{definition} \cite{bowditch}
A pair $(G, \mathcal{H})$ is relatively hyperbolic if $G$ admits a geometrically finite action on a proper, hyperbolic space $X$ such that the set $\mathcal{H}$ consists of exactly the maximal parabolic subgroups and each of these are finitely generated.
\end{definition}

We are now ready to state the result we are interested in; how well specified the geometry of these spaces are, given the group and peripheral group structure.  What we find is that while the core of the space is well-defined up to quasi-isometry, the shape of the cusps can break quasi-isometry between candidate spaces.  Here `core' means the space that is the lift of a connected compactum that separates the ends in the quotient, The inspiration for this result came from observations in \cite{sisgroman}, which asserted a version of Lemma \ref{key} in Remark A.6.

\begin{thmz}
Any relatively hyperbolic group with infinite peripheral subgroups acts as above on hyperbolic spaces that are not equivariantly quasi-isometric.

\label{equivariance}
\end{thmz}

We need one more definition before stating the heavy-lifting lemma.

\defi \cite{groman}
For a connected, locally finite metric graph $\Gamma$ with edge lengths 1, and increasing function $f: \mathbb{R}_{\geq 0} \rightarrow \mathbb{R}_{\geq 0}$ that is coarsely at least exponential, the associated \emph{combinatorial horoball} $\mathcal{C}_f(\Gamma)$ is a graph with vertex set 
$$ V(\cah(\Gamma)) := \Gamma^0 \times \bbn$$ where the points $(v,n)$ and $(v,n+1)$ are connected by edges of length 1, and each level $ \Gamma^0 \times {n}$ has an edge of length 1 between them if their distance in $\Gamma$ was less than or equal to $f(n)$.
\edefi

This object is mostly used for groups, in which case the combinatorial horoball of a (sub)group $G$ will be denoted $\mathcal{C}_f(G)$ and refer to $\mathcal{C}_f(\Gamma(G))$ for some understood Cayley graph $\Gamma$.  In the case of a subgroup, it will be assumed the intended graph is the natural subgraph of $\Gamma(G)$.

\begin{lemma}
Let $f(x) := 2^x$ and $g(x) := 2^{2^x}$.  Then for any finitely generated infinite group $H$, the combinatorial horoballs $\mathcal{C}_f(H)$ and $\mathcal{C}_g(H)$ are not quasi-isometric. 
\label{key}
\end{lemma}

\pf
Here we might be tempted to apply the idea, explained for example in \cite{bh}, that the growth of balls in a graph is a quasi-isometry invariant, because these two graphs by design have different growth rates.  However, this statement is made specifically for graphs which arise as Cayley graphs of finitely generated groups.  Implicitly in this formulation, we are using the assumption that our graph has uniformly bounded valence, which regrettably is not true for these combinatorial horoballs.  We must do a little more work.  We will denote distance in $\mathcal{C}_f(H)$ by $d_f$, distance in $\mathcal{C}_g(H)$ by $d_g$, and distance restricted to the zero level of the horoball (which is independent of which scaling function is used) by $d_H$.

The first thing we observe about these spaces is that they are $\delta-$hyperbolic for some $\delta$, by \cite{sisgroman} Appendix A, and the boundaries are single points.  In the geodesic ray definition of the boundary, these points are the equivalence class of rays that point straight `upwards', consisting entirely of vertical edges. Therefore any $(c,c)$ quasi-isometry $\phi: \mathcal{C}_f(H) \rightarrow \mathcal{C}_g(H)$ between the combinatorial horoballs, which acts by homeomorphism on the boundaries of hyperbolic spaces, must take these geodesic rays to quasi-geodesic rays in the equivalence class of the lone boundary element on the right which are in turn $B$ close to geodesic rays for some value $B$.  Up to bounded distance, then, we can assume such a map takes a ray $\{ (x,i) \ | \ i \geq 0 \}$ to some geodesic ray $\{ (y,i) \ | \ i \geq k \}$, for $x,y \in H$.  Because of this potential error between geodesic and quasi-geodesic rays, the next paragraph is performed up to a bounded constant $B=B(\delta)$ which can be added to the appropriate constant found for each application of the quasi-isometry.

We proceed by contradiction.  Let $\phi$ be a $(c,c)$ quasi-isometry between these spaces, and denote its quasi-inverse by $\psi$.  Without loss of generality, assume $\psi$ also has $(c,c)$ for quasi-isometry constants.  We first need to note that the zero level of $\mathcal{C}_f(H)$ has bounded height in the image.  To do this, it is sufficient to consider what happens to the vertices.  Pick an arbitrary point in this level set, $(x,0)$.  By quasi-inverses, $\psi \circ \phi ((x,0))$ has bounded distance from $(x,0)$, and so has height bounded by some multiple of $c$.  Due to the above observation of geodesic rays, $\{ (x,i) \ | \ i \geq 0 \}$ goes to some geodesic ray $\{ (y,i) \ | \ i \geq k \}$ and under the quasi inverse $\psi$,  $\{ (y,i) \ | \ i \geq k \}$ goes to some geodesic ray $\{ (z,i) \ | \ i \geq \ell \}$.  In particular, this means that the value of $\ell$ is \emph{at least} $\frac{k}{c} - {c}$.  Because the height of $\psi \circ \phi ((x,0))$ is bounded above by some multiple of $c$, this says that the value of $k$ is also bounded above by some function of $c$.  So the height of any point $\phi ((x,0))$ is bounded by a function of $c, \delta$.  A symmetric argument guarantees the same is true for $\psi((y,0))$.  Call the maximum of these height bounds $D$.

Now consider two points $x_0, x_N \in H$ with horospherical distance $d_H (x_0,x_N) = N$.  We can pick these for any value of $N$ desired by the assumption that $H$ has infinite diameter.  Subdivide an $H-$geodesic between these points into a path $x_0, x_1 \ldots x_N$ so that each successive point is at distance 1; note that each $x_i$ will necessarily be a vertex.  Let $(y_i, h_i) := \phi ((x_i,0))$, recalling that $h_i < D$.  Because we know that 
$$ d_g(\phi((x_i,0)),\phi((x_{i+1},0))) \leq 1c + c = 2c  $$
and
$$ d_g(\phi((x_i,0)),(y_i,0)) \leq D + B $$
the triangle inequality guarantees that 
$$ d_g((y_i,0),(y_{i+1},0)) \leq 2c + 2(D+B) $$

Recall that the last term comes from the discrepancy between geodesics and quasi-geodesics, and only depends on the hyperbolicity constant.  This observation implies that $d_H(y_i,y_{i+1})$ is uniformly bounded (independent of $i, N$) by some constant we label $E$.  So we know that for any choice of $x_0,x_N$, the distance $d_H(y_0,y_N) \leq EN$.  Then we also know by the way we defined $\mathcal{C}_g(H)$ that

\begin{align*}
d_H(y_0,y_N) \leq EN \implies & d_g((y_0,0),(y_N,0)) \leq 2 \lceil \log_2 (\log_2(EN)) \rceil + 3 \\
\implies & d_g(\phi((x_0,0)),\phi((x_N,0))) \leq 2 \lceil \log_2 (\log_2(EN)) \rceil + 3 +2B \quad (\dagger)
\end{align*}

This is because we can adapt Lemma 3.10 from \cite{groman} to observe that geodesics between vertices in these combinatorial horoball spaces will always consist of traveling towards the boundary point along vertical edges, traveling along at most 3 horizontal edges, and traveling again vertically downwards, which achieves at most the distance listed above.  However, the fact that $\phi$ is a $(c,c)$ quasi-isometry dictates that

\begin{align*}
 d_g(\phi((x_0,0)),\phi((x_N,0))) \geq & \frac{1}{c} d_f((x_0,0),(x_N,0)) -c \\
\geq &  \frac{1}{c} (2 \lfloor \log_2(N)\rfloor +1) -c \quad \quad (\dagger \dagger)
\end{align*}

by another application of 3.10 from \cite{groman}.  We see that for sufficiently large values of $N$, the statements $\dagger$ and $\dagger \dagger$ are incompatible, contradicting the existence of such a map $\phi$.

\epf

\begin{myproof} \textbf{of \ref{equivariance}} \\
Let $(G,H)$ be a relatively hyperbolic pair, such that $H$ is infinite.  We construct two spaces, $X_1, X_2$ as follows.  $X_i$ will be a copy of $\Gamma(G)$, the Cayley graph, with combinatorial horoballs glued on to all cosets of $H$.  In $X_1$, allow the scaling function to be $2^n$, and in $X_2$ allow the scaling function to be $2^{2^n}$.  Again by \cite{sisgroman}, we note that this resultant space is hyperbolic for both cases.  The equivalence of definitions of relative hyperbolicity tell us that these spaces are acted upon in the appropriate sense of \cite{bowditch}. \\

In order to apply our lemma in the correct way, we need to show that any equivariant map between the $X_i$ will coarsely take cusps to cusps.  Denote by $Q_1, Q_2$ the quotient of $X_1, X_2$ respectively by the action of $G$.  By the assumption that these spaces are hyperbolic, they have a well defined boundary, and by the assumption that the action is geometrically finite, these will be spaces that have finitely many, isolated boundary points, in 1-1 correspondence with the number of (conjugacy classes of) peripheral subgroups.  \\

Suppose $f$ is a G-equivariant (c,c) quasi-isometry $f : X_1 \rightarrow X_2$.  Then $f$ descends to a quasi-isometry $f_q: Q_1 \rightarrow Q_2$.  Explicitly, let $p_i : X_i \rightarrow Q_i$ be the quotient maps.  Then we can define $f_q (q_1) = p_2 ( f( p_1^{-1} (q_1) ) )$.  Note this is well-defined because of the assumption of equivariance of the map $f$.  We claim $f_q$ is a $(c,c)-$quasi-isometry.  Let $z_1, z_2$ be points in $Q_1$.

\begin{align*}
d_{Q_2}( f_q(z_1), f_q(z_2) ) &= d_{Q_2}( p_2( f( p_1^{-1}(z_1))), p_2( f( p_1^{-1}(z_2) ))) & \text{By definition of } f_q \\
&\leq  d_{X_2}(  f( p_1^{-1}(z_1)),  f( p_1^{-1}(z_2) )) & \text{Projection does not} \\
& & \text{increase distance} \\
&\leq c \ d_{X_1}( p_1^{-1}(z_1),  p_1^{-1}(z_2) ) + c& f \text{ is a (c,c)-QI} \\
&=  c \ d_{Q_1}(z_1,  z_2) + c & \text{with careful choice of pre-image}
\end{align*}
\\
and \\
\begin{align*}
d_{Q_2}( f_q(z_1), f_q(z_2) ) &= d_{Q_2}( p_2( f( p_1^{-1}(z_1))), p_2( f( p_1^{-1}(z_2) ))) & \text{By definition of } f_q \\
&=  d_{X_2}(  f( p_1^{-1}(z_1)),  f( p_1^{-1}(z_2) )) & \text{Choosing pre-image points} \\
& & \text{at the minimum distance} \\
&\geq \frac{1}{c}d_{X_1}( p_1^{-1}(z_1),  p_1^{-1}(z_2) ) - c& f \text{ is a (c,c)-QI} \\
&\geq  \frac{1}{c}d_{Q_1}(z_1,  z_2) - c  & \text{Projection does not} \\
& & \text{increase distance} 
\end{align*}

Finally, to satisfy the quasi-onto condition, we note that any point in $X_2$ is bounded distance from the image of $f$, so the same will be true (with the same bound) when the points are projected downstairs. \\

Now, because $f_q$ is a QI from $Q_1$ to $Q_2$, it must act as a homeomorphism on the discrete boundary.  In particular, it takes geodesic rays representing these boundary points, to infinite length quasi-geodesic rays in $Q_2$.  The only such rays that exist in this space are those that represent the cusps that are the quotient of peripheral groups upstairs.  Therefore, downstairs, $f$ coarsely maps  cusps to cusps.  Because of how we defined this map, this forces $f$ to map, again coarsely, our combinatorial horoballs in $X_1$ to those in $X_2$.  Then, by Lemma \ref{key}, this map cannot be a QI after all, so we have a contradiction, meaning no such map can exist.
\end{myproof}

It is conjectured by the author that we may drop equivariance in the statement of \ref{equivariance} if we allow the scaling functions to be super-exponential and super-super-exponential, with the proof of Lemma \ref{key} being similar just with more 2s.  In this case, we would expect to find that we naturally cannot coarsely map cusps into the `core' of the target space by a divergence argument.
\end{section}

\begin{section}{Decay of Rigidity}
These two results tell us that as we loosen the conditions that we use to classify negatively curved groups, we also lose some of the metric structure and end behavior their corresponding spaces enjoy.  We sum up this meta-statement in the following table, where `Yes' indicates that structure is rigid

\begin{table}[h]

\begin{center}
\begin{tabular}{ | c | c || c | c | }
\hline
 Group Property & Action Type on & Boundary/ & QI Type \\[1ex]
 & Hyperbolic $X$ & Limit Set & of X \\
\hline \hline
 Hyperbolic & Geometric & Yes & Yes \\ [1ex]
Relatively Hyperbolic & Geometrically Finite & Yes &  No \\[1ex]
Acylindrically Hyperbolic & Universal & No &  No \\
\hline
\end{tabular}
\end{center}
\caption[Summary of Rigidity]{Quasi-isometric and Limit Set Rigidity of Hyperbolicity Generalizations}
\label{rigidity}
\end{table}
It should be noted here that in the third column, we are referring to the well-defined Gromov boundary of the group in the first row and the Gromov boundary of the `cusped space' in the second row, which is well-defined by \cite{bowditch}.
\end{section}

\bibliography{refs}
\bibliographystyle{alpha}
\end{document}